\documentclass{article}

\usepackage{template_arxiv/arxiv}

\usepackage[utf8]{inputenc} 
\usepackage[T1]{fontenc}    
\usepackage[hidelinks]{hyperref}       
\usepackage{url}            
\usepackage{booktabs}       
\usepackage{amsfonts}       
\usepackage{nicefrac}       
\usepackage{microtype}      
\usepackage{lipsum}		
\usepackage{graphicx}
\usepackage{natbib}
\usepackage{doi}

\usepackage{amsmath,amsthm,amssymb,bm}
\usepackage{xcolor}
\usepackage{subfig}
\usepackage{appendix}
\usepackage{algorithm}
\usepackage{algpseudocode}
\usepackage{adjustbox}
\usepackage[normalem]{ulem}
\usepackage{xparse}

\newcommand{\norm}[1]{\left\lVert#1\right\rVert}

\makeatletter
\NewDocumentCommand{\LeftComment}{s m}{%
  \Statex \IfBooleanF{#1}{\hspace*{\ALG@thistlm}}\(\triangleright\) #2}
\makeatother

	\title{On the implementation of Adaptive and Filtered MHE}

	\author{ \href{https://orcid.org/0000-0003-4694-6339}{\includegraphics[scale=0.06]{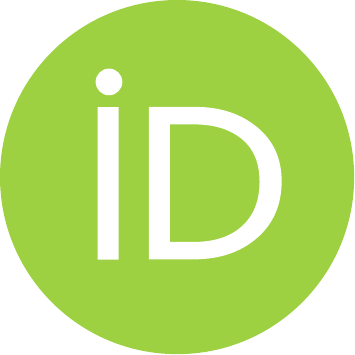}\hspace{1mm}Federico Oliva} \\
	Dipartimento di Ing. Civile e Ing. Informatica \\
	University of Rome ``Tor Vergata''
	\texttt{lvofrc95@outlook.it} \\
	\And
	{Daniele Carnevale} \\
	Dipartimento di Ing. Civile e Ing. Informatica \\
	University of Rome ``Tor Vergata''
	\texttt{carnevaledaniele@gmail.com} \\
}

	\renewcommand{\headeright}{Technical Report}

	\hypersetup{
		pdftitle={On the implementation of Adaptive and Filtered MHE},
		pdfsubject={Optimisation and Control},
		pdfauthor={Federico Oliva, Daniele Carnevale},
		pdfkeywords={Adaptive observer, Nonlinear systems, Moving Horizon Estimator},
	}

\begin{document}
\maketitle

\begin{abstract}
	Optimisation-based algorithms known as \textit{Moving Horizon Estimator} (MHE) have been developed through the years. This paper illustrates the implementation of the policy introduced in the companion paper submitted to the 18th IFAC Workshop on Control Applications of Optimization \citep{Oliva01}, in which we propose two techniques to reduce the computational cost of MHEs. These solutions mainly rely on output filtering and adaptive sampling. The use of filters reduces the total amount of data used by MHE, shortening the length of the moving window (buffer) and consequently decreasing the time consumption for plant dynamics integration. Meanwhile, the proposed adaptive sampling policy discards those sampled data that do not allow a sensible improvement of the estimation error. Algorithms and numerical simulations are provided to show the effectiveness of the proposed strategies.
\end{abstract}



\section{Introduction}

Observer design represents a well-known task in control theory, playing a pivotal role in the controller design. A classic example of the interconnections of these two problems is described by the \textit{Separation Principle} stated by Kalman in its most famous and elegant definition. As far as observers are concerned, many different approaches have been studied and developed over the years; from the most used linear algorithms like the Luenberger and Kalman observers \citep{Davis,thrun}, up to more involved solutions for nonlinear systems \citep{thrun,Reif,Carnevale,Karagiannis}. Another topic of great interest in the field of observers is the adaptive approach, whose task consists in estimating model parameters too \citep{Luders,Krener,Marino,Tyukin,Marino2,Marino3}. A fascinating approach addressing a general observation problem (both adaptive and non-adaptive) is the one proposed by the Moving Horizon Estimators (MHE) \citep{Michalska,Grizzle,Kang2006,Wynn2014,Sui2010,Suwantong2014,Schiller2021}. These are very powerful optimisation-based observers, but they suffer from computational cost issues. These issues make their online implementation quite hard to reach. To increase the computational speed, we propose in  \citep{Oliva01} to exploit output-measurement filters to reduce the number of model integrations needed to solve the optimisation problem at the basis of standard MHE algorithms. Moreover, we propose an adaptive-sampling policy, making the observer capable of autonomously choosing the "best" measurements to be considered in the moving window buffer. The strategies proposed hinge upon MHE 
convergence analysis presented in \citep{Aeyels,Kang2006,Glad,Menini_2019,Menini_2022}. 

\medskip
This short work describes in detail the implementation of solutions proposed in \citep{Oliva01}, presenting the Algorithms to implement the \textit{Filtered} MHE and the \textit{Adaptive} MHE and providing as well some simulation results. More specifically, in \ref{scn:convergence_analysis} the \textit{Standard} MHE frameworks is described, while \ref{scn:strategies} treats both the \textit{Filtered} MHE and \textit{Adaptive} MHE cases. Conclusions are drawn in \ref{scn:conclusions}, as well as future developments.


\section{Moving Horizon Estimators}
\label{scn:convergence_analysis}

Consider the general framework described in \citep{Oliva01}, namely a nonlinear system in the following form

\begin{subequations}
	\label{eqn:conv_general_system}
	\begin{align}
		\dot{\bm{\xi}} & = f(\bm{\xi},\bm{u})  \\
		\bm{y}         & = h(\bm{\xi},\bm{u}).
	\end{align}
\end{subequations}

\medskip
In this set of equations $\bm{\xi}\in\mathbb{R}^n$ is the state vector, $\bm{u}\in \mathbb{R}^m$ is the control input, and $\bm{y} \in \mathbb{R}^p$ is the measured output. From this general framework, in order to define the MHE, system \eqref{eqn:conv_general_system} shall be considered in its \textit{N-lifted} form. Accordingly to \citep*{Grizzle} and \citep{Tousain}, and following the notation introduced in \citep{Oliva01}, the \textit{N-lifted} system generated from \eqref{eqn:conv_general_system} is described by an output operator in the following form:

\begin{align}
	\text{H}_k(\bm{\xi}_{k-N+1},u(\cdot)) \!=\!\!\!
	\begin{bmatrix}
		\! h(\bm{\xi}_{k-N+1},\bm{u}_{k-N+1})                                         \\
		\! h(\phi(t_{k-N+2},t_{k-N+1},\bm{\xi}_{k-N+1},\bm{u}(\cdot)),\bm{u}_{k-N+2}) \\
		\! \vdots                                                                     \\
		\! h(\phi(t_{k},t_{k-1},\bm{\xi}_{k-1},\bm{u}(\cdot)),\bm{u}_{k})
	\end{bmatrix},
	\label{eq:Hmap}
\end{align}

\medskip
where $(\text{Y}_k,\text{U}_k)$ are the output and input sample buffers of length $N$ and down-sampling $N_{T_s}$. Moreover, $\phi(t_k;t_{k-1},\bm{\xi}_{k-1},\bm{u}(\tau)) = \bm{\xi}(t_k) = \bm{\xi}_{k}$ describes the solution of \eqref{eqn:conv_general_system} at time $t_k$ with $\tau\in [t_{k-1},t_k]$, and $\bm{u}(t)$ is assumed to be a piece-wise constant function of the time, within sampling times of length $T_s$. Therefore, The MHE associated to \eqref{eqn:conv_general_system} through its \textit{N-lifted} system can be defined as

\begin{subequations}
	\begin{align}
		\label{eqn:conv_genral_observer}
		\bm{{\zeta}}_{i+1} & = \bm{\Psi}(\text{Y}_k,\text{U}_k,\text{H}_k(\bm{{\zeta}}_i,\bm{u}(\cdot)),
		\bm{\zeta}_i),                                                                                   \\
		\bm{\zeta}_{0}     & =\hat{\bm{\xi}}_{k-N+1},\text{ with, } i=0\dots K, \ K\in\mathbb{N},        \\
		\hat{\bm{\xi}}_{k} & = \phi(t_k,t_{k-N+1},\bm{\zeta}_K,\bm{u}(\cdot)),\label{eqn:propagation}
	\end{align}
\end{subequations}

\medskip
In this set of equations we consider $\hat{\bm{\xi}}_k$ as the estimated state vector of $\bm{\xi}_k$. Instead, $i$ is the iteration number of the optimisation algorithm defining the updated value of $\bm{\zeta}$. This value is computed through $\Psi$, that describes a general algorithm designed to solve optimisation problems (e.g. simplex or gradient-based solutions). This optimisation algorithm is iteratively launched within the intervals $[t_{k-1},t_{k}]$ for $K$ times. To sum up, the estimation problem at the basis of MHEs consists in finding the solution to the following minimisation problem:

\begin{equation}
	\label{eqn:conv_system_optimisationProblem}
	\underset{\hat{\bm{\xi}}_{k-N+1}}{\text{min}} \text{V}_k(\text{Y}_k,\text{U}_k,\hat{\text{H}}_k,\hat{\bm{\xi}}_{k-N+1}).
\end{equation}

\medskip
The cost $V_k$ has been defined as a quadratic function, namely

\begin{equation}
	\label{eqn:conv_V_general}
	V_k(\text{Y}_k,\text{U}_k,{\text{H}}_k,\zeta) \triangleq \sum\limits_{j=1}^{N} (\text{Y}_k^j - \hat{\text{H}}_k^j)^TW_j(\text{Y}_k^j - \hat{\text{H}}_k^j),
\end{equation}

\medskip
where we considered the notation $\hat{H}_k = \text{H}_k(\hat{\bm{\xi}}_{k-N+1},$ $\bm{u}(\cdot))$. Moreover, $W_i\in\mathbb{R}^{p\times p}$ are symmetric and positive definite weight matrices and $(\text{Y}_k^j,\hat{\text{H}}_k^j)$ are the \textit{j-th} rows of the matrices $Y_k$ and $\hat H_k$, respectively. For a detailed discussion on the \textit{Standard} MHE structure of the optimisation problem solution and on its convergence properties, please refer to \citep{Oliva01}. In addition to these considerations it is important to remark that  if the plant dynamics \eqref{eqn:conv_genral_observer} are linear, then it is possible to consider the Newton algorithm yielding zero estimation error with just $K=2$ since $V_k$ is quadratic, and therefore resulting in the \textit{dead-beat} observer described in \citep*{Grizzle}\footnote{The discrete LTI plant of \eqref{eqn:conv_general_system} need to be evaluated selecting $N_{T_s}=1$.}.

\medskip
Lastly, with reference to the \textit{Standard} MHE proposed in \citep{Oliva01}, a more clear description of the algorithm is proposed in pseudo-code in \ref*{alg:MHE_standard}. Remember that $K$ represents the number of optimisation steps performed by the optimisation algorithm $\bm{\Psi}$ chosen to solve \eqref{eqn:conv_system_optimisationProblem}. Moreover, note that the integral performed in the last algorithm step on the system model is computed numerically, as the analytical solution $\phi$ is not in general available in closed-form.

\begin{algorithm}
	\caption{Standard MHE}\label{alg:MHE_standard}
	\begin{algorithmic}
		\For {$t_k=k\cdot N_{T_s}\cdot T_s \in [t_0, t_f]$}
		\State \LeftComment{Update buffers}
		\State $\text{Y}_k[1:N-1] \gets \text{Y}_k[2:N]$,
		$\text{U}_k[1:N-1] \gets \text{U}_k[2:N]$
		\State $\text{Y}_k[N] \gets y_k$, $\text{U}_k[N] \gets u_k$
		\State \LeftComment{Initialization of
			\eqref{eqn:conv_genral_observer}}
		\State $\zeta_0 \gets \hat{\bm{\xi}}_{k-N+1}$
		\State \LeftComment{Optimisation procedure}
		\For{$i=0,\dots,K-1$}
		\State {\small $\zeta_{i+1} \gets \bm{\Psi}(\text{Y}_k,\text{U}_k,\hat{\text{H}}_k,\zeta_i)$}
		\EndFor
		\State $\hat{\bm{\xi}}_{k-N+1} \gets \zeta_K$
		\State \LeftComment{Propagate the estimate evaluating \eqref{eqn:propagation}}
		\For{$i\leq N-1$}
		\State {\small $\hat{\bm{\xi}}_{k-N+1+i} \gets \int\limits_{t_{k-N+i}}^{t_{k-N+1+i}}f(\hat{\bm{\xi}}(\tau),u(\tau))d\tau$}
		\EndFor
		\EndFor
		\State \LeftComment{Estimate at
			$t= t_{k + q T_s} \neq t_k$, with $q\in\mathbb{Z}$,
			$|q|< N_{T_s}$}
		\State $\hat{\bm{\xi}}(t) \gets \int\limits_{t_{k}}^{t_{k}+qT_s}
			f(\hat{\bm{\xi}}(\tau),u(\tau))d\tau$
	\end{algorithmic}
\end{algorithm}


\section{Strategies for computational efficiency}
\label{scn:strategies}

This section considers the \textit{Filtered} MHE and \textit{Adaptive} MHE described in \citep{Oliva01}, and aims at better describing the code implementation. Moreover, some considerations regarding the estimation convergence are reported. Generally speaking the focus is on the role of parameters $N$ (buffer length) and $N_{T_s}$ (down-sampling). Some considerations will be explained by simulation results on two different models, the Van der Pol oscillator and a nonlinear system describing the dynamics of plasma waves in presence of runaway electrons \citep{Buratti,Carnevale}.

\subsection*{Filtered MHE}
\label{subscn:filtered_considerations}

As reported in \citep{Oliva01} a first solution to decrease the computational cost necessary to solve \eqref{eqn:conv_system_optimisationProblem} consists in adding a filtered version of the system output measurements $y_k$. Roughly speaking, by increasing the conditions available at each time instants, the buffer length $N$ can be reduced. Results in this sense have been presented on a Van der Pol oscillator.


For the model setup and main results refer to \citep{Oliva01}. The major effect of adding a filtering action on the output $y$ is that the filter state vector shall be propagated, updated and stored too during the optimisation process. Indeed, a general filter can be described as a discrete-time system, namely

\begin{subequations}
	\label{eqn:filter_generic_DSS}
	\begin{align}
		\bm{\xi_f}(k+1) & = \gamma(\xi_f(k),y(k)),  \\
		y_f(k+1)        & = \beta(\bm{\xi_f}(k+1)),
	\end{align}
\end{subequations}

where $\bm{\xi_f} \in \mathbb{R}^l$ is the filter state, $y\in\mathbb{R}$ is the signal measurement to be filtered, and $y_f\in\mathbb{R}$ is the actual filtered signal. The pseudo-code relative to the \textit{Filtered} MHE is presented in \ref*{alg:MHE_filters}. The notation is the same used for \textit{Standard} MHE.

\begin{algorithm}
	\caption{Filtered MHE}\label{alg:MHE_filters}
	\begin{algorithmic}
		\For {$t_k=k\cdot N_{T_s}\cdot T_s \in [t_0, t_f]$}

		\State \LeftComment{Output filter}
		\State $\bm{\xi_f}(k) \gets \gamma(\bm{\xi_f}(k-1),y_k)$
		\State $y_k^f \gets \beta(\bm{\xi_f}(k))$
		\State $\overline{y}_k \gets [y_k \ y_k^f]^T$

		\State \LeftComment{Update buffers}
		\State $\text{Y}_k[1:N-1] \gets \text{Y}_k[2:N]$,
		$\text{U}_k[1:N-1] \gets \text{U}_k[2:N]$
		\State $\text{Y}_k[N] \gets \overline{y}_k$, $\text{U}_k[N] \gets u_k$
		\State \LeftComment{Initialization of
			\eqref{eqn:conv_genral_observer}}
		\State $\zeta_0 \gets \hat{\bm{\xi}}_{k-N+1}$
		\State \LeftComment{Optimisation procedure}
		\For{$i=0,\dots,K-1$}
		\State {\small $\zeta_{i+1} \gets \mathbf{\Psi}(\text{Y}_k,\text{U}_k,\hat{\text{H}}_k,\zeta_i)$}
		\EndFor
		\State $\hat{\bm{\xi}}_{k-N+1} \gets \zeta_K$
		\State \LeftComment{Propagate the estimate evaluating \eqref{eqn:propagation}}
		\For{$i\leq N-1$}
		\State {\small $\hat{\bm{\xi}}_{k-N+1+i} \gets \int\limits_{t_{k-N+i}}^{t_{k-N+1+i}}f(\hat{\bm{\xi}}(\tau),u(\tau))d\tau$}
		\State $\hat{\bm{\xi}}_{\bm{f}}( k-N+1+i) \gets \sum\limits_{i=k-N+i}^{k-N+1+i}\gamma(\hat{\bm{\xi}}_{\bm{f}}(i),\hat{y}_i)$
		\EndFor
		\EndFor
		\State \LeftComment{Estimate at
			$t= t_{k + q T_s} \neq t_k$, with $q\in\mathbb{Z}$,
			$|q|< N_{T_s}$}
		\State $\hat{\bm{\xi}}(t) \gets \int\limits_{t_{k}}^{t_{k}+qT_s}
			f(\hat{\bm{\xi}}(\tau),u(\tau))d\tau$
		\State $\hat{\bm{\xi}}_{\bm{f}}(k+q) \gets \sum\limits_{i=k}^{k+q}\gamma(\hat{\bm{\xi}}_{\bm{f}}(i),\hat{y}_i)$
	\end{algorithmic}
\end{algorithm}

\subsection*{Adaptive MHE}
\label{subscn:adaptive_MHE}

The second solution proposed in \citep{Oliva01} to speed up MHEs implementation consists in defining a policy according to which the observer can automatically select the output measurements that are most informative for the solution of \eqref{eqn:conv_system_optimisationProblem}. This selection is done accordingly to the following indices $\delta_k$ and $d_{\text{V}}$, that are a good representative of the output signal richness as well as the precision of the estimate:

\begin{subequations}
	\label{eqn:down_PE}
	\begin{align}
		\Sigma_k     & = \{\sigma_i \triangleq \norm{\text{Y}_k^{i+1} - \text{Y}_k^i} \ \text{ s.t. } i\in [1, \dots,N-1]\}, \\
		\delta_k     & = \sum\limits_{i=1}^{N-1}\sigma_i + \norm{y_k-\text{Y}_k^1},                                          \\
		d_{\text{V}} & = \norm{\text{Y}_k - \hat{\text{H}}_k}.
	\end{align}
\end{subequations}

\medskip
The adaptive policy is implemented by thresholding both $\delta_k$ and $d_{\text{V}}$ and deciding whether to run or not the estimation accordingly. Moreover, if no estimation is performed for a long time due to a high precision reached, the output buffer is re-initialised to avoid excessively time-consuming model integrations (Parameter $N_{\text{max}}$ \citep{Oliva01}). Again, the pseudo-code relative to the \textit{Adaptive} MHE is presented in \ref*{alg:MHE_adaptive}. Results have been presented on the following model, describing the dynamics of plasma waves amplitude ($\xi_1$) and anisotropy of the runaway electrons velocity distribution ($\xi_2$) given by \citep{Buratti}:

\begin{subequations}
	\label{eqn:runaway_model}
	\begin{align}
		\dot{\xi}_1 & = \epsilon(-2\xi_1\xi_2-2S+Q),                                                              \\
		\dot{\xi}_2 & = \epsilon(-\nu \xi_2 + \xi_3(\xi_1\xi_2 + S)-\gamma_1\dfrac{\xi_2}{1+\dfrac{\xi_2}{W_t}}), \\
		\dot{\xi_3} & = 0,                                                                                        \\
		y           & = \xi_2,
	\end{align}
\end{subequations}

\begin{figure}[h!]
	\centering
	\includegraphics[height=6cm,keepaspectratio]{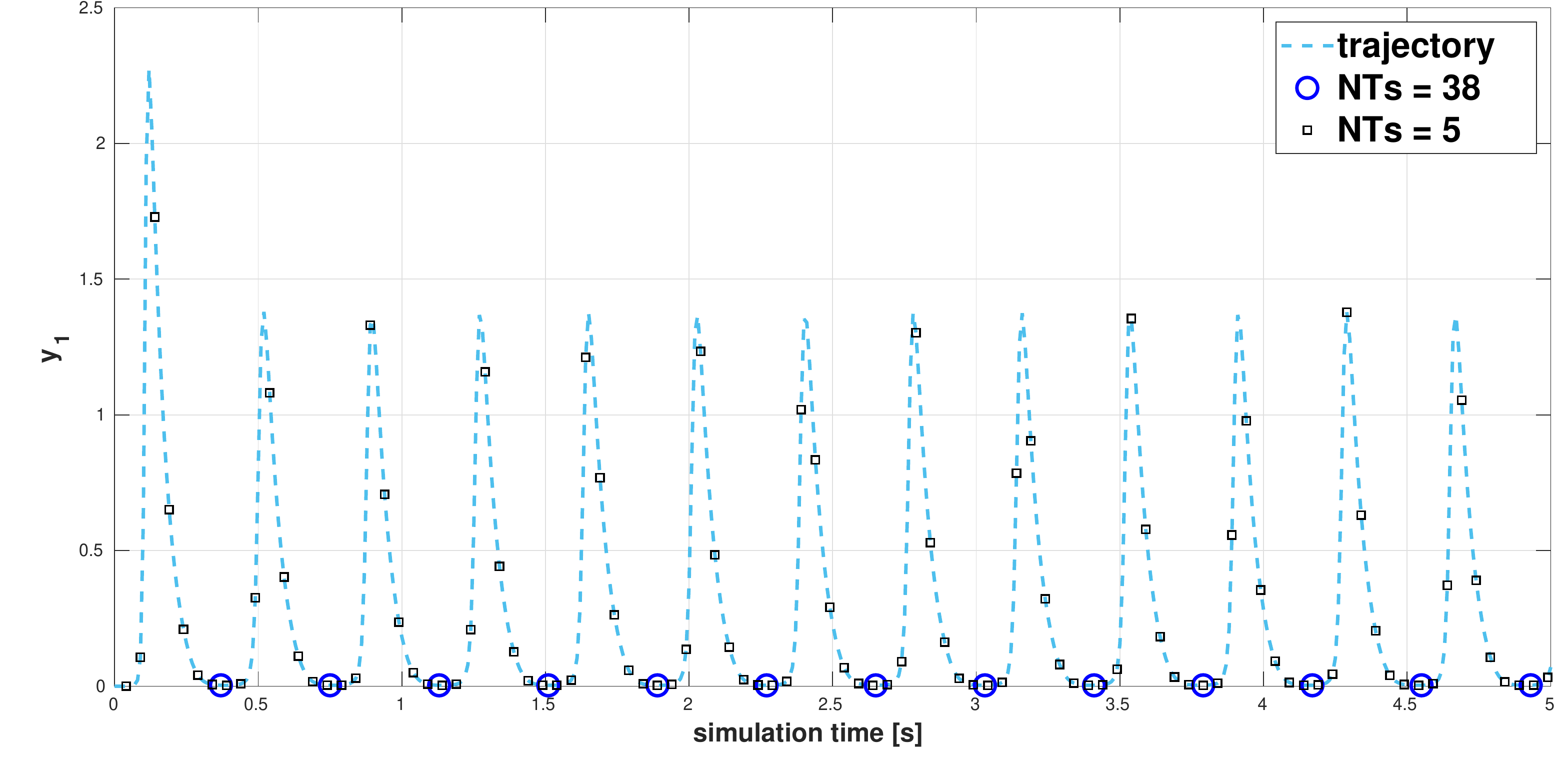}
	\caption{This figure compares the system trajectory (dashed line) and the sampled measurements when $N_{T_s} = 5$ (squares) and $N_{T_s} = 38$ (circles).}
	\label{fig:down_traj_Nts}
\end{figure}

\medskip
where $y \in \mathbb{R}$ is the output, $\bm{\xi}\in\mathbb{R}^3$ is the state vector. \autoref{fig:down_traj_Nts} shows the actual measurement samples, when different values of fixed $N_{T_s}$ are considered, namely $N_{T_s} = 5$ and $N_{T_s} = 38$. As reported in \citep{Oliva01}, when $N_{T_s} = 38$, index $\delta_k$ is low and the state is not correctly estimated. Indeed, this means that the sampled trajectory is not highly informative. This can be clearly seen in \autoref{fig:down_traj_Nts} where the blue circles have a nearly null value compared to the general trend of the trajectory. Clearly, an informative set of measurements shall be sampled over the trajectory peaks, rather than on the flat and dense intervals between them.

\begin{figure}[h!]
	\centering
	\includegraphics[height=6cm,keepaspectratio]{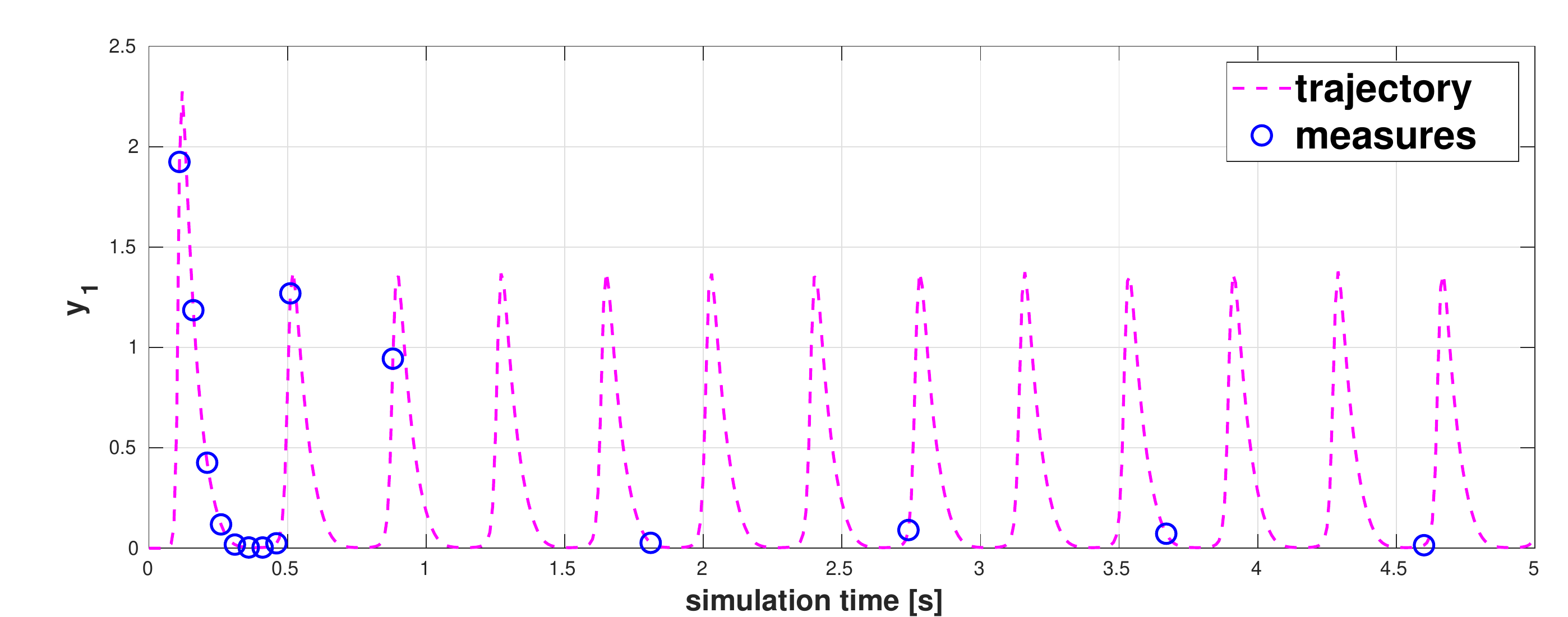}
	\caption{This figure shows the true state trajectory for system \eqref{eqn:runaway_model} (black dashed line), and the sampled measures provided to the observer when the adaptive sampling policy is used (blue big circles).}
	\label{fig:down_adaptive_sampling}
\end{figure}

The results of the actual \textit{Adaptive} MHE are presented in \autoref{fig:down_adaptive_sampling}. The simulation follows the very same setup considered in \citep{Oliva01}. Note how the sampling is now concentrated at the beginning of the simulation and on the peaks, in order to quickly reduce the estimation error and then correct the estimation when the measurements were informative enough. Whenever the estimation error increases, the observer samples the output in correspondence of those trajectory regions in which the signal richness is higher. The estimation precision is comparable to the \textit{Standard} MHE, yet significantly reducing the computational cost.

\begin{algorithm}
	\caption{Adaptive MHE}\label{alg:MHE_adaptive}
	\begin{algorithmic}

		\For {$t_k=k\cdot N_{T_s}\cdot T_s \in [t_0, t_f]$}

		\LeftComment{Adaptive sampling}
		\If{AC1 $\And$ AC2}

		\State \LeftComment{Check max spacing}
		\If{$T(y_k) - T(\text{Y}_k^N) \geq N_{\max}$}
		\State \LeftComment{Re-initialise buffer}
		\For{$i\leq N$}
		\State $\text{Y}_k^i = y_{k-N+i}$
		\EndFor
		\Else
		\State \LeftComment{Update buffers}
		\State $\text{Y}_k[1:N-1] \gets \text{Y}_k[2:N]$
		\State $\text{U}_k[1:N-1] \gets \text{U}_k[2:N]$
		\State $\text{Y}_k[N] \gets y_k$, $\text{U}_k[N] \gets u_k$
		\EndIf

		\State \LeftComment{Initialization of
			\eqref{eqn:conv_genral_observer}}
		\State $\zeta_0 \gets \hat{\bm{\xi}}_{k-N+1}$
		\State \LeftComment{Optimisation procedure}

		\For{$i=0,\dots,K-1$}
		\State {\small $\zeta_{i+1} \gets \mathbf{\Psi}(\text{Y}_k,\text{U}_k,\hat{\text{H}}_k,\zeta_i)$}
		\EndFor
		\State $\hat{\bm{\xi}}_{k-N+1} \gets \zeta_K$
		\State \LeftComment{Propagate the estimate evaluating \eqref{eqn:propagation}}

		\For{$i\leq N-1$}
		\State {\small $\hat{\bm{\xi}}_{k-N+1+i} \gets \int\limits_{t_{k-N+i}}^{t_{k-N+1+i}}f(\hat{\bm{\xi}}(\tau),u(\tau))d\tau$}
		\EndFor
		\EndIf
		\EndFor
		\State \LeftComment{Estimate at
			$t= t_{k + q T_s} \neq t_k$, with $q\in\mathbb{Z}$,
			$|q|< N_{T_s}$}
		\State $\hat{\bm{\xi}}(t) \gets \int\limits_{t_{k}}^{t_{k}+q T_s}
			f(\hat{\bm{\xi}}(\tau),u(\tau))d\tau$

	\end{algorithmic}
\end{algorithm}

\subsection*{Filtered and Adaptive MHE }
\label{subscn:jont_MHE}

This last section focuses on the joint action of \textit{Filtered} MHE and \textit{Adaptive} MHE. In fact, \textit{Filtered} MHE reduces the computational speed of each optimisation solution of \eqref{eqn:conv_system_optimisationProblem}, while \textit{Adaptive} MHE reduces the number of optimisation required to reach a good estimation precision. Therefore, if both \ref{alg:MHE_filters} and \ref{alg:MHE_adaptive} were used in combination, namely exploiting both output filtering and adaptive sampling, we would both decrease the single optimisation time and the total number of performed optimisations. Indeed, the speed up factor would be even greater.

\begin{figure}[h!]
	\centering
	\includegraphics[height=6cm,keepaspectratio]{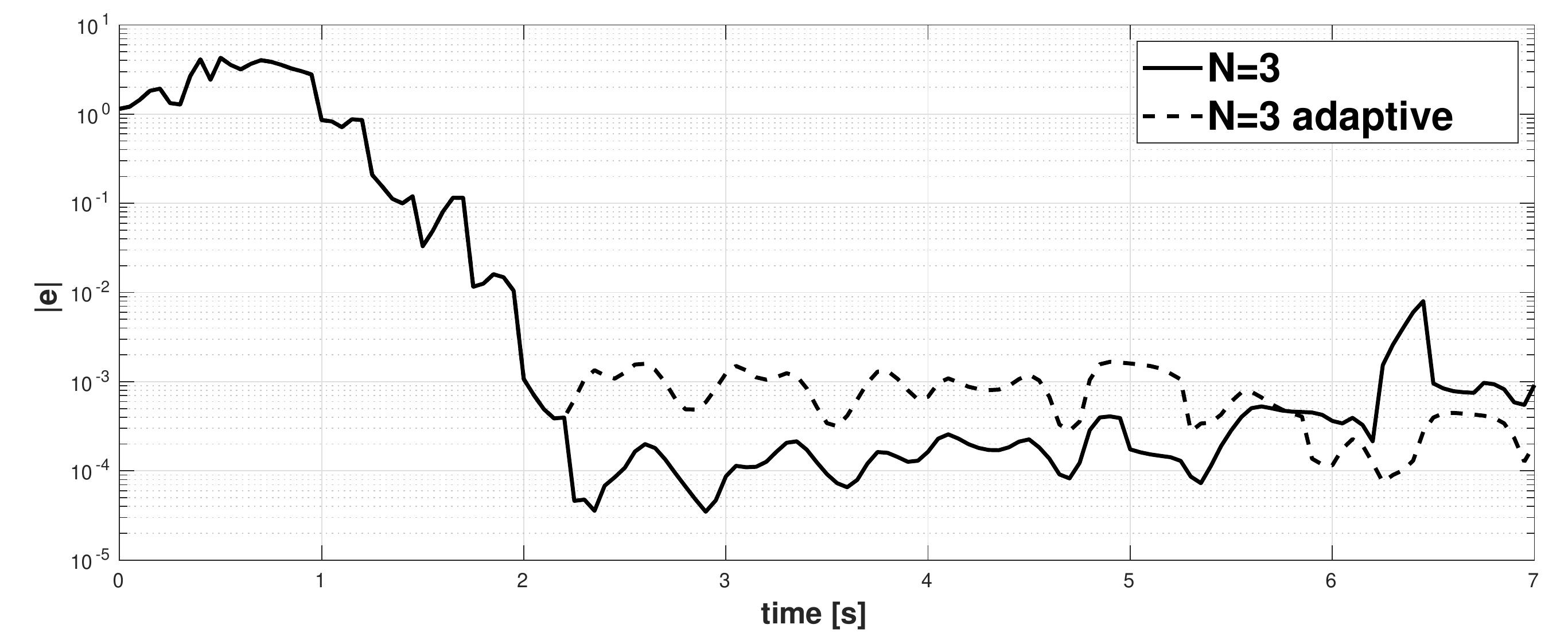}
	\caption{This figure shows the estimation error on both the filtered (solid line) and filtered-adaptive sampling (dashed line) cases on system \eqref{eqn:double_pendulum} estimation.}
	\label{fig:combined_error}
\end{figure}

\medskip
In order to test the performance of both output filtering and adaptive sampling we combined the \textit{Filtered} MHE and the \textit{Adaptive} MHE on a double pendulum, considering $K=40$ and two output filter, namely a dirty-derivative and an integrator with loss [cfr. Oliva et al. 2022]. The model considered for this set of simulations consists of a double pendulum. The general structure of the equations resembles the usual model used in robotics, namely

\begin{equation}
	\label{eqn:double_pendulum}
	\mathbb{M}(q)\Ddot{q} + \mathbb{V}(q,\dot{q})\dot{q} + \mathbb{G}(q) = \tau,
\end{equation}

where $q = (\theta_1, \theta_2)$ are the angular positions of the links, the $\mathbb{M}$ term describes the system inertia, the $\mathbb{V}(q,\dot{q})$ term the friction and Coriolis effects, and $\mathbb{G}(q)$ the gravitational force. Each pendulum was considered with length $L = 1m$, and mass $M=1kg$. The system state is $\bm{\xi} = [q_1, q_2, \dot{q}_1, \dot{q}_2] \in \mathbb{R}^4$, while the output measurement is $y = q_1 \in \mathbb{R}$. The estimation error norm is presented in logarithmic scale in \autoref{fig:combined_error}.

The registered total computation time was $25.4$s for the joint action of \textit{Filtered} MHE and \textit{Adaptive} MHE which is nearly $60\%$ lower than the simple \textit{Filtered MHE} case reported in \citep{Oliva01}.


\section{Conclusions and future work}
\label{scn:conclusions}

This work builds on the results presented in \citep{Oliva01}. In the first section \ref*{scn:convergence_analysis}, the \textit{Standard} MHE structure was recalled and presented more accurately in pseudo-code layout. Moreover, some remarks on the convergence properties introduced in \citep{Oliva01} are reported. More specifically, the special case of an MHE turning into a \textit{dead-beat} observer has been considered. In \ref*{scn:strategies} the solutions proposed in \citep{Oliva01} to speed up MHE are considered. More specifically, the importance of filter propagation is stressed, and the \textit{Filtered} MHE is presented in pseudo-code. As far as \textit{Adaptive} MHE is concerned, the simulation results on plasma dynamics model described in \citep{Oliva01} are consolidated with some considerations on both the fixed and adaptive sampling policies. Lastly, the \textit{Adaptive} MHE is presented in pseudo-code as well.

\medskip
Future developments of this tool will consider a detailed analysis of the effect of measurement noise and model uncertainties on the observer performance and selection of filters. Another interesting point to be developed is using this tool for control laws design, relating it to the general MPC solution. Indeed, output tracking problems could be solved by including control parameters in the augmented state. For instance, classic MRAC or state feedback frameworks could be reproduced. Clearly, this application shall be considered as separated from the state estimation problem, because no separation principle is available for nonlinear systems.

\subsubsection*{Acknowledgemts}
Thanks to Corrado Possieri and Mario Sassano for their comments and fundamental insights into the development of this work.

\bibliographystyle{unsrtnat}
\bibliography{main}

\end{document}